\documentclass[11pt,a4paper]{article}
\usepackage[latin1]{inputenc}
\usepackage{amsfonts,amsmath,amssymb,theorem,mathrsfs}
\usepackage{graphicx,epic,color}
\usepackage{bbm}

\usepackage{hyperref}

\newtheorem{theorem}{Theorem}
\newtheorem{proposition}{Proposition}
\newtheorem{lemma}{Lemma}
{\theorembodyfont{\rmfamily} \newtheorem{remark}{Remark}}

\newcommand{\N}{\mathbb{N}}
\newcommand{\R}{\mathbb{R}}
\newcommand{\bbone}{\mathbbm{1}}
\newcommand{\bbeta}{\boldsymbol{\beta}}

\title{A population model for $\Lambda$-coalescents with neutral mutations}
\author{Andreas Nordvall Lagerås\thanks{Stockholm University, Department of Mathematics, division of Mathematical Statistics, 10691 Stockholm, Sweden, e-mail: \href{mailto:andreas@math.su.se}{andreas@math.su.se}}}

\begin{document}

\maketitle

\begin{abstract}
Bertoin and Le Gall \cite{BertoinLeGall:2003} introduced a certain probability measure valued Markov process that describes the evolution of a population, such that a sample from this population would exhibit a genealogy given by the so-called $\Lambda$-coalescent, or coalescent with multiple collisions, introduced independently by Pitman \cite{Pitman:1999} and Sagitov \cite{Sagitov:1999}. We show how this process can be extended to the case where lineages can experience mutations. Regenerative compositions enter naturally into this model, which is somewhat surprising, considering a negative result by Möhle \cite{Mohle:2005b}.
\end{abstract}


\section{Introduction}

A coalescent with multiple collisions, or $\Lambda$-coalescent, $\Pi=(\Pi_t)_{t\geq 0}$, is a Markov process on the space $\mathscr{P}(\N)$, the partitions of $\N=\{1,2,\dots\}$, such that for all $n$, $\Pi^{(n)}$, the restriction of $\Pi$ to $[n]=\{1,\dots,n\}$, is a Markov process with the following transitions: If $\Pi^{(n)}_t$ has $b$ blocks, then any collection of $k$ blocks, coalesce into one block at rate $\lambda_{b,k}$ for $2\leq k \leq b \leq n$. Note that the rate only depends on the number of blocks, not their sizes. By considering $\Pi^{(n)}$ and $\Pi^{(n+1)}$ one realizes that $\lambda_{b,k}=\lambda_{b+1,k}+\lambda_{b+1,k+1}$. From this it follows, see \cite{Pitman:1999}, that
$$
\lambda_{b,k}=\int_{[0,1]}x^{k-2}(1-x)^{b-k}\Lambda(dx)\text{ for }2\leq k \leq b,
$$
for some finite measure $\Lambda$ on $[0,1]$.

$\Lambda$-coalescents were introduced independently by Pitman \cite{Pitman:1999} and Sagitov \cite{Sagitov:1999} as a generalization of the Kingman coalescent \cite{Kingman:1982}. All these processes can arise as limiting processes when studying the genealogy for a finite sample of individuals from a haploid (one parent per child) population, see \cite[Proposition 7]{Schweinsberg:2000}.

Consider a large population with constant size $N$ for all generations, which are furthermore non-overlapping. A sample of $n$ individuals, labeled by $[n]$, form a partition by grouping together those who have had a common ancestor by generation $t$ backwards in time, i.e.\ those whose lineages have coalesced into a common lineage by that time. The $\Lambda$-coalescent, restricted to $[n]$, is a possible limiting process when $N\to\infty$, and time and the distribution of the number of children of each individual are scaled properly. Furthermore, it is possible to obtain a coalescent with simultaneous multiple collisions, but we will in this paper not consider such so-called $\Xi$-coalescents, see \cite{MohleSagitov:2001, Schweinsberg:2000} for more details.

The Kingman coalescent is a $\Lambda$-coalescent with $\Lambda=\delta_0$, i.e.\ with the only type of transition being a merger of two blocks at a time. This process is the natural limiting process for many population models, roughly speaking those models where the number of children for each individual always is small compared to the total population size as the size tends to infinity. The probability of more than two lineages in the sample coalescing at the same time is then negligible in the limit. A $\Lambda$-coalescent with $\Lambda\neq\delta_0$ corresponds to a population where occasionally single individuals have offspring constituting a positive fraction of the entire next generation as the population size tends to infinity. If several of the lineages in your sample belong to that fraction, they will coalesce into a single lineage at that moment.

The main result of this paper is a description of the dynamics of the \emph{whole} population when all lineages experience neutral mutations, i.e.\ mutations that do not influence the chance of survival. Earlier studies have only described how introducing mutations influences the dynamics of the genealogy of a \emph{sample} from the population.

We will proceed as follows. First, in Section \ref{section:paintboxes_and_bridges}, we will acquaint ourselves with useful representations of random partitions and coalescent processes. Here we will also find a description of a population model such that the genealogy of a sample from this population is given by the $\Lambda$-coalescent. When we introduce mutations in the population, an obvious way of partitioning a sample of individuals is by their common genotypes. In Section \ref{section:mutations_in_the_sample}, a general recursion formula is given for the distribution of the family sizes in the sample. Section \ref{section:regenerative_composition_structures} might at first be seen as a detour into the theory of regenerative composition structures, i.e.\ a special kind of ordered partitions, especially since it is known that our type of partitions can never appear from these regenerative composition structures if one simply disregards the order in the composition. This theory, however, is used in the last Section \ref{section:mutations_in_the_population}, in which we present a model for the whole population, and not just a sample from it, when all lineages experience neutral mutations, such that the distribution of a sample from this population is in accordance with the result in Section  \ref{section:mutations_in_the_sample}.

\section{Paintboxes and a population process}\label{section:paintboxes_and_bridges}

In this section we will see how one can use probability measures to construct random partitions of $\N$.  Let $p$ be a probability measure with atoms of sizes $\mathbf{b}=(b_i)_{i\in\N}$ in non-increasing order. Let $(R_i)_{i\in\N}$ be an i.i.d.\ sample from $p$ and define the equivalence relation $\sim_p$ on $\N$ by $i\sim_pj$ if $R_i=R_j$. Then the equivalence classes of $\sim_p$ form the sought partition of $\N$.

A more common way of using such a sequence $\mathbf{b}$, is to partition $[0,1]$ into intervals $(I_i)_{i\in\N_0}$ with lengths $(b_i)_{i\in\N_0}$, where $b_0:=1-\sum_ib_i$, and let the equivalence relation $\sim_{\mathbf{b}}$ on $\N$ be defined by $i\sim_{\mathbf{b}}j$ if $V_i,V_j\in I_k$ for some $k\geq 1$, where $(V_i)_{i\in\N}$ are i.i.d.\ $U(0,1)$, and let $i$ such that $V_i\in I_0$ be in equivalence classes of their own. In general, one can use a random measure $\pi$, and carry out the construction pointwise, given $\pi=p$. This is equivalent to using atoms of random sizes $\bbeta=(\beta_i)_{i\in\N}$, and the construction is called a paintbox construction, see \cite{Kingman:1982}.

A random partition of $\N$ is called exchangeable if its restriction to $[n]$ has a distribution that is invariant under permutations of the labels $[n]$ for all $n\in\N$. For example, if $(\Pi_t)_{t\geq 0}$ is a $\Lambda$-coalescent, then $\Pi_t$ is exchangeable for all $t$. Kingman has shown, see e.g.\ \cite{Kingman:1982}, that any exchangeable random partition of $\N$ can be obtained from a paintbox construction, e.g.\ with $\bbeta=(\beta_i)_{i\in\N}$ being the almost sure limit of $(l_i(n)/n)_{i\in\N}$, where $l_i(n)$ is the size of the $i$th largest block of the partition restricted to $[n]$. Here, and elsewhere in this paper, we understand limits to be taken as $n\to\infty$, unless otherwise indicated.

If one enumerates the blocks of $\Pi_s=\{A^s_1,A^s_2,\dots\}$, and for $t>s$ considers the blocks of $\Pi_t=\{A^t_1,A^t_2,\dots\}$, then each $A^t_i=\cup_{j\in C^{s,t}_i}A^s_j$ for some $C_i^{s,t}$. It is a property of coalescent processes that the partition $\Pi^{s,t}=\{C^{s,t}_1,C^{s,t}_2,\dots\}$ is exchangeable and distributed as $\Pi_{t-s}$. Bertoin and Le Gall \cite{BertoinLeGall:2003} showed that there exists a collection $(\pi_{s,t})_{s\leq t}$ of random probability measures on $[0,1]$ such that $\sim_{\pi_{s,t}}$ corresponds to $\Pi^{s,t}$ for all $s<t$. For fixed $s$ and increasing $t$ this collection describes the genealogy of the population further and further backwards in time. The random partition $\Pi^{s,t}$ should be interpreted as describing how the lineages present in the population at time $s$ coalesce into lineages at time $t$, and it thus has a meaning even with negative arguments $s<t$, since this corresponds to future events in the population relative to time zero. We shall study, and later extend, the dynamics of the Markov process $\rho=(\rho_t)_{t\geq 0}:=(\pi_{-t,0})_{t\geq 0}$, which describes the evolution of the population forwards in time. Heuristically, $\rho_t(dr)$ represents the descendants at time $t$ of the fraction $dr$ of the population at time zero.

If $\Lambda(0)=0$, then the dynamics of $\rho$ can be described by a measure $\nu(dx):=\Lambda(dx)/x^2$. Let $\{(\tau_i,X_i,U_i)\}_{i\in\N}$ be a Poisson process on $\R\times(0,1]\times(0,1)$ with intensity measure $dt\otimes\nu(dx)\otimes du$. Assume for the moment $|\nu|:=\nu((0,1])<\infty$. We will use this Poisson point process to construct the process $\rho$. We let $(\tau_i)_{i\in\N}$, which by the assumption is a homogeneous Poisson process with rate $|\nu|$, be the jump times of $\rho$. At a time $\tau_i$, the conditional law of $\rho_{\tau_i}$, given $\rho_{\tau_i-}$, is
\begin{equation}\label{rho_jump}
(1-X_i)\rho_{\tau_i-}+X_i\delta_{R_i},
\end{equation}
where $X_i$ has distribution  $\nu/|\nu|$ and $R_i$ is a sample from $\rho_{\tau_i-}$, picked by the inverse transformation method: $R_i:=\inf(r:\rho_{\tau_i-}([0,r])>U_i)$. Between jumps, $\rho$ remains constant.

The heuristic interpretation of these dynamics is that a person is chosen from the population just before the jump at time $\tau_i$, and she is identified  with her family labeled $R_i$. At the time of the jump, she begets offspring of proportion $X_i$ of the total population, and we say that the \textit{litter} $i$, born at $\tau_i$, has size $X_i$. The rest of the population must thus be scaled down by a factor $1-X_i$ and the atom corresponding to her family is increased with mass $X_i$.

Bertoin and Le Gall \cite[Corollary 2]{BertoinLeGall:2003} showed the following.
\begin{proposition}\label{BertoinLeGall} If $(\nu_n)_{n\in\N}$ is a sequence of finite measures, with $\Lambda_n(dx)=x^2\nu_n(dx)$ converging weakly to a finite measure $\Lambda$ on $[0,1]$, then the sequence of processes $(\rho^n)_{n\in\N}$, where each $\rho^n$ is governed by the respective $\nu_n$, converges in distribution, in the sense of weak convergence of finite-dimensional marginals, to the process $\rho$ corresponding to the collection $(\pi_{s,t})_{s\leq t}$ associated with the limiting $\Lambda$-coalescent. If $\Lambda(0)=\Lambda(1)=0$ and $\int_{(0,1)}x\nu(dx)<\infty$, then the convergence can be strengthened to almost sure convergence.
\end{proposition}
Thus, $\rho$ has a meaning even in the case $|\nu|=\infty$, and in particular, the description of $\rho$ above, for the case $|\nu|<\infty$, can be extended to the case $\int_{(0,1)}x\nu(dx)<\infty$.

\section{Mutations in the sample}\label{section:mutations_in_the_sample}

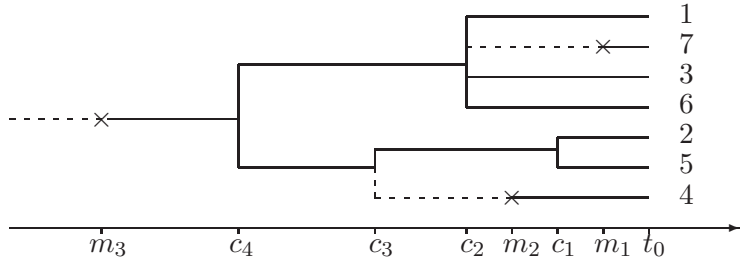
\begin{figure}
\centering
\setlength{\unitlength}{0.8cm}
\begin{picture}(12,5)
\put(0,1){\vector(1,0){12}}
\dashline{0.1}(0,2.8)(1.5,2.8)
\dashline{0.1}(6,1.5)(8.25,1.5)
\dashline{0.1}(7.5,4)(9.75,4)
\dashline{0.1}(6,1.5)(6,2)
\drawline(1.6,2.8)(3.75,2.8)
\drawline(3.75,2)(6,2)
\drawline(3.75,3.7)(7.5,3.7)
\drawline(6,2.3)(9,2.3)
\drawline(7.5,3)(10.5,3)
\drawline(7.5,3.5)(10.5,3.5)
\drawline(7.5,4.5)(10.5,4.5)
\drawline(8.25,1.5)(10.5,1.5)
\drawline(9,2)(10.5,2)
\drawline(9,2.5)(10.5,2.5)
\drawline(9.75,4)(10.5,4)
\drawline(3.75,2)(3.75,3.7)
\drawline(6,2)(6,2.3)
\drawline(7.5,3)(7.5,4.5)
\drawline(9,2)(9,2.5)
\put(1.5,2.8){\makebox(0,0){$\times$}}
\put(8.25,1.5){\makebox(0,0){$\times$}}
\put(9.75,4){\makebox(0,0){$\times$}}
\drawline(1.5,1)(1.5,0.9)
\put(1.3,0.6){$m_3$}
\drawline(3.75,1)(3.75,0.9)
\put(3.6,0.6){$c_4$}
\drawline(6,1)(6,0.9)
\put(5.9,0.6){$c_3$}
\drawline(7.5,1)(7.5,0.9)
\put(7.4,0.6){$c_2$}
\drawline(8.25,1)(8.25,0.9)
\put(8.1,0.6){$m_2$}
\drawline(9,1)(9,0.9)
\put(8.9,0.6){$c_1$}
\drawline(9.75,1)(9.75,0.9)
\put(9.6,0.6){$m_1$}
\drawline(10.5,1)(10.5,0.9)
\put(10.4,0.6){$t_0$}
\put(11,4.4){1}
\put(11,3.9){7}
\put(11,3.4){3}
\put(11,2.9){6}
\put(11,2.4){2}
\put(11,1.9){5}
\put(11,1.4){4}
\end{picture}
\caption{A family tree for a sample of 7 individuals. $t_0$ is present time. Chronological time increases to the right, whereas the time of the coalescent increases to the left. Coalescence of lineages occurs at times $c_1,\dots,c_4$. Mutations are denoted by $\times$ and occur at times $m_1,m_2$ and $m_3$. The partition of this sample into families is $\big\{\{1,2,3,5,6\},\{4\},\{7\}\big\}$.}\label{FamilyTree}
\end{figure}

In a population genetics setting, it is natural to introduce \textit{mutations} along the lineages and ask how the individuals of your sample are partitioned into different families according to their genotype. We assume that mutations always give rise to new types of individuals never seen before in the population (the so-called infinite alleles model), and that when tracing the lineages backwards in time, there is a constant intensity $\mu$ per lineage for a mutation to occur, i.e.\ if we draw the family tree of the sample, the mutations constitute a homogeneous Poisson process with intensity $\mu$ along each branch, see Figure \ref{FamilyTree} for an example.

A quantity of interest is $q(\mathbf{a})$,  $\mathbf{a}=(a_1,a_2,\dots)$, the  probability of observing a partition with $a_i$ families of size $i$. When we trace the genealogy of the lineages backwards in time, the probability of a mutation to occur first is $\mu n/(\mu n + \lambda_n)$, and the probability of a collision of $k$ lineages happening first is $\binom{n}{k}\lambda_{n,k}/(\mu n + \lambda_n)$, for $k=2,\dots,n$. By the Markov property of the $\Lambda$-coalescent, we can condition on the type of event that happens first, and obtain a recursion for $q(\mathbf{a})$. Möhle \cite{Mohle:2005a} was the first to provide this recursion for $\Lambda$- (and even $\Xi$-) coalescents. Let $\mathbf{e}_k$ be the $k$th unit vector in $\R^{\infty}$ and $\lambda_n:=\sum_{k=2}^n\lambda_{n,k}$.

\begin{proposition}\label{mohle} {\rm({\bf Möhle's recursion} \cite{Mohle:2005a})} $q(\mathbf{e}_1) = 1 $, and 
\begin{equation}\label{Mohles_recursion}
q(\mathbf{a}) = \frac{\mu n}{\mu n + \lambda_n}q(\mathbf{a}-\mathbf{e}_1)+\sum_{k=2}^n\frac{\binom{n}{k}\lambda_{n,k}}{\mu n + \lambda_n}\sum_{j=1}^{n-k+1}\frac{j(a_j+1)}{n-k+1}q(\mathbf{a}+\mathbf{e}_j-\mathbf{e}_{j+k-1}),
\end{equation}
for $n=\sum_i ia_i >1$, where $q(\mathbf{a})=0$ if any $a_i<0$.
\end{proposition}
There are no known closed formulas solving \eqref{Mohles_recursion} for general $\Lambda$, except for the cases $\Lambda=\delta_0$ (Kingman's coalescent) and $\Lambda=\delta_1$.

The parts of this formula should be interpreted as follows. If a mutation occurs first, the rest of the sample is described by $\mathbf{a}-\mathbf{e}_1$. If a merger of $k$ lineages occurs first, and it occurs in a family represented by $j+k-1$ lineages, then after that merger, the sample will consist of $n-k+1$ lineages and be described by $\mathbf{a}+\mathbf{e}_j-\mathbf{e}_{j+k-1}$. In particular, there will now be $a_j+1$ families of size $j$. The probability that the merger of $k$ lineages affected a family of size $j+k-1$ is given by $j(a_j+1)/(n-k+1)$, since the merger could have resulted in any of the $j$ lineages in any of the $a_j+1$ families with equal probability by the exchangeability. We refer to Möhle \cite{Mohle:2005a} and Dong et al.\ \cite{DongEtAl:2006} for more detailed discussions. The latter paper extends coalescent processes, so that their blocks become \emph{frozen} when they encounter a mutation, and then do not partake in the further evolution of the process. With frozen blocks enclosed by $\langle$ and $\rangle$, the path of this process, realized as in Figure \ref{FamilyTree}, would be
\begin{align*}
&\big\{\{1\},\dots,\{7\}\big\}\overset{m_1}{\to} \big\{\{1\},\dots,\{6\},\langle7\rangle\big\}\overset{c_1}{\to} \big\{\{1\},\{2,5\},\{3\},\{4\},\{6\},\langle 7\rangle\big\}\\
&\overset{m_2}{\to} \big\{\{1\},\{2,5\},\{3\},\langle 4\rangle,\{6\},\langle 7\rangle\big\}\overset{c_2}{\to} \big\{\{1,3,6\},\{2,5\},\langle 4\rangle,\langle 7\rangle\big\}\\
&\overset{c_4}{\to} \big\{\{1,2,3,5,6\},\langle 4\rangle,\langle 7\rangle\big\}\overset{m_3}{\to} \big\{\langle 1,2,3,5,6\rangle,\langle 4\rangle,\langle 7\rangle\big\}.
\end{align*}
The partition into families is obtained when all blocks are frozen.

\section{Regenerative composition structures}\label{section:regenerative_composition_structures}

Most results in this section are from Gnedin and Pitman \cite{GnedinPitman:2005a}. A partition of $n\in\N$ is an unordered collection of natural numbers $\{n_1,\dots,n_k\}$ such that $n_1+\cdots+n_k=n$. An ordered partition is called a \emph{composition}, and we say that $n_1,\dots,n_k$ are its parts. A composition structure $\mathscr{C}$ is a sequence $(C_n)_{n\in\N}$ of random compositions of $n$ such that if $n$ balls are distributed into an ordered series of boxes according to $C_n$, then $C_{n-1}$ is obtained by discarding one of the balls picked uniformly at random, and deleting an empty box in case one is created. A composition structure is \emph{regenerative} if for all $n\geq m\geq 1$, given that the first part is $m$, the remaining composition of $n-m$ is distributed as $C_{n-m}$.

\begin{figure}
\centering
\setlength{\unitlength}{1cm}
\begin{picture}(11,2.2)(-1,-1.1)

\thicklines
\drawline(0,-0.1)(0,0.1)
\drawline(10,-0.1)(10,0.1)
\put(-1,-0.1){$S$}
\put(-0.1,0.2){0}
\put(9.9,0.2){1}

\drawline(0,0)(1,0)
\drawline(1.5,0)(2.5,0)
\drawline(3.5,0)(4.5,0)
\drawline(5.5,0)(6,0)
\drawline(7,0)(7.5,0)
\drawline(8,0)(8.5,0)
\drawline(8.8,0)(9.3,0)
\dottedline{0.1}(9.3,0)(9.9,0)

\drawline(0,-1.1)(0,-0.9)
\drawline(10,-1.1)(10,-0.9)
\put(-1,-1.1){$S^{(z)}$}
\drawline(0,-1)(1.1,-1)
\drawline(3.3,-1)(4.4,-1)
\drawline(5.6,-1)(6.7,-1)
\drawline(7.3,-1)(8.4,-1)
\dottedline{0.1}(8.4,-1)(9.9,-1)
\thinlines

\dottedline{0.1}(0,-1)(5.5,0)
\dottedline{0.1}(1.1,-1)(6,0)
\dottedline{0.1}(3.3,-1)(7,0)
\dottedline{0.1}(4.4,-1)(7.5,0)
\dottedline{0.1}(5.6,-1)(8,0)
\dottedline{0.1}(6.7,-1)(8.5,0)
\dottedline{0.1}(7.3,-1)(8.8,0)
\dottedline{0.1}(8.4,-1)(9.3,0)
\dottedline{0.1}(10,-0.9)(10,-0.1)

\put(1.2,0.2){\vector(0,-1){0.2}}
\put(1.1,0.3){$V_4$}
\put(2,0.2){\vector(0,-1){0.2}}
\put(1.9,0.3){$V_7$}
\put(2.7,0.2){\vector(0,-1){0.2}}
\put(2.6,0.3){$V_2$}
\put(3.3,0.2){\vector(0,-1){0.2}}
\put(3.2,0.3){$V_5$}
\put(4.8,0.2){\vector(0,-1){0.2}}
\put(4.7,0.3){$z$}
\put(6.15,0.2){\vector(0,-1){0.2}}
\put(6.00,0.3){$V_1$}
\put(6.55,0.2){\vector(0,-1){0.2}}
\put(6.40,0.3){$V_3$}
\put(6.9,0.2){\vector(0,-1){0.2}}
\put(6.85,0.3){$V_6$}
\end{picture}
\caption{An illustration of sampling with $V_1,\dots,V_7$ from $S$ resulting in $(n_1,n_2,n_3,n_4)=(1,1,2,3)$, and how to construct $S^{(z)}$ from $S$.}\label{Rz}
\end{figure}
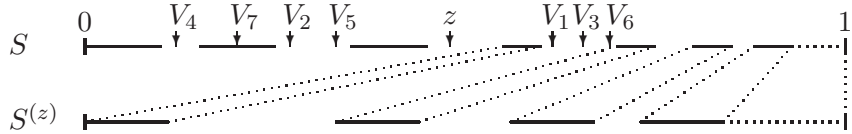

We will see that one can obtain regenerative compositions with the appropriate sampling procedure. Let $(V_i)_{i\in\N}$ be i.i.d.\ $U(0,1)$ and let $(V_{in})_{i\in[n]}$ be the ordered sample of $(V_i)_{i\in[n]}$, meaning $V_{1n}\leq \cdots\leq V_{nn}$. Given a closed set $S\subseteq [0,1]$, we can construct a composition $C_n$ as follows. Partition $[n]$ into blocks of \emph{consecutive} integers by letting $j$ and $j+1$ be in different blocks if $[V_{jn},V_{j+1,n}]\cap S\neq\emptyset$. Let the parts of $C_n$ be given by the sizes of the blocks in increasing order of their elements, see Figure \ref{Rz}. We will in general also allow a random closed set $\mathcal{S}\subseteq [0,1]$, independent of $(V_i)_{i\in\N}$, where the construction is carried out given the realization $\mathcal{S}=S$. We then say that $C_n$ is obtained by sampling from $\mathcal{S}$.

For any closed set $S\subseteq [0,1]$ and $z\in[0,1)$, define $D(S,z):=\inf S\cap(z,1]$, where we let $\inf\emptyset:=1$. For $S$ and $z$ such that $D(S,z)<1$, define
\begin{equation*}
S^{(z)}:=\left\{\frac{x-D(S,z)}{1-D(S,z)}:x\in S\cap[D(S,z),1]\right\}
\end{equation*}
This is the part of $S$ strictly to the right of $D(S,z)$ scaled back to $[0,1]$, see Figure \ref{Rz}. We say that a random closed set $\mathcal{S}\subseteq [0,1]$ is \emph{multiplicatively regenerative} if for each $z\in[0,1)$, given $D(\mathcal{S},z)<1$, the set $\mathcal{S}^{(z)}$ is independent of $[0,D(\mathcal{S},z)]\cap\mathcal{S}$, and has the same distribution as $\mathcal{S}$.

Let $\{(\tau_i,X_i)\}_{i\in\N}$ be a Poisson process on $\R_+\times (0,1)$ with intensity measure $dt\otimes\nu(dx)$ for a measure $\nu$ with $\int_{(0,1)}x\nu(dx)<\infty$, and let $\mu\geq 0$ be a constant. The notation here is intensionally similar to the one in the previous sections of this paper, but we assume for the moment no relation to these. We call the process $Z=(Z_t)_{t\geq 0}$ a \emph{multiplicative subordinator} with characteristics $(\mu,\nu)$ if
$$
Z_t:=1-e^{-\mu t}\prod_{i:\tau_i\leq t}(1-X_i),
$$
for all $t\geq 0$. The name is justified by the property that $(1-Z_{t'})/(1-Z_t)$ has the same distribution as $1-Z_{t'-t}$ and is independent of $(Z_u)_{0\leq u\leq t}$ for $t'>t$. We obtain an ordinary subordinator by the transformation $Z_t\mapsto -\log(1-Z_t)$.\begin{remark}\label{mom_cond}
We need the moment condition on $\nu$ above to obtain a non-trivial process $Z$, since if $\int_{(0,1)}x\nu(dx)=\infty$, then $\sum_{i:\tau_i\leq t}\log(1-X_i)=-\infty$ almost surely for $t>0$, see Campbell's Theorem in \cite[p.\ 28]{Kingman:1993}, and thus $\prod_{i:\tau_i\leq t}(1-X_i) = 0$ for all $t>0$.
\end{remark}

Let $\mathcal{R}$ be the closed range of the multiplicative subordinator $Z$. Proposition \ref{exp_sampling} collects some results of \cite{GnedinPitman:2005a}.
\begin{proposition}\label{exp_sampling}
The closed range $\mathcal{R}$ of $Z$ is  multiplicatively regenerative, and conversely, all multiplicatively regenerative sets can be seen as the range of some multiplicative subordinator, whose characteristics are determined up to a positive constant. Sampling from $\mathcal{R}$ produces a regenerative composition structure $\mathscr{C}$, and all regenerative composition structures can be obtained by sampling from a regenerative set.
\end{proposition}
Since we have these relations between regenerative composition structures, multiplicative subordinators, and multiplicatively regenerative sets, we  also say that $(\mu,\nu)$ are the characteristics of the regenerative composition structure $\mathscr{C}$ of the proposition. In particular, the probability of the first part having size $m$ in $C_n$, is $q(n:m)=\Phi(n:m)/\Phi(n)$, where
\begin{equation}\label{qnm}
\Phi(n:m)=\mu n1(m=1)+\binom{n}{m}\int_0^1x^m(1-x)^{n-m}\nu(dx),
\end{equation}
and $\Phi(n)=\sum_{m=1}^n\Phi(n:m)$. We see that the characteristics $(\mu,\nu)$ and $(c\mu,c\nu), c>0$, produce the same regenerative composition structure. We will need more detailed results about the first part of a regenerative composition $C_n, n\geq 2$. It can have size one if either $V_{1n}\in\mathcal{R}$, or $V_{1n}\notin\mathcal{R}$ and $\mathcal{R}\cap[V_{1n},V_{2n}]\neq\emptyset$. The expressions for the following probabilities are taken from the proof of Theorem 5.2 on p. 457 of \cite{GnedinPitman:2005a}.
\begin{align}
q(n:1)'&:=P(V_{1n}\in\mathcal{R})=\frac{\mu n}{\Phi(n)},\label{qn1prim}\\
q(n:1)''&:=P(V_{1n}\notin\mathcal{R},[V_{1n},V_{2n}]\cap\mathcal{R}\neq\emptyset)=\frac{n}{\Phi(n)}\int_0^1x(1-x)^{n-1}\nu(dx).\label{qn1bis}
\end{align}
Möhle \cite[Theorem 3.1]{Mohle:2005b} showed 
\begin{proposition}\label{mohle_neg}
The recursion \eqref{Mohles_recursion} cannot be solved by a partition obtained by disregarding the order of the parts of a regenerative composition structure, unless $\Lambda$ has all its mass in either 0 or 1.
\end{proposition}

\section{Mutations in the population}\label{section:mutations_in_the_population}

For a population without mutations, $\rho_t\to\delta_{\mathbf{e}}$ in distribution as $t\to\infty$, where $\mathbf{e}$ has distribution $U[0,1]$ and is called the primitive Eve \cite[Proposition 1 and Definition 4]{BertoinLeGall:2003}, so that all of the population belongs to the primitive Eve's family. This is a sort of genetic drift where, by chance, some genotype eventually makes up the whole population. When mutations are possible, no such absorbing state exists since new mutations appear, and we can hope for the existence of a non-trivial stationary distribution of $\rho$.

We shall now investigate what happens with $\rho$, describing the evolution of the population forwards in time, when individual lineages mutate at constant rate $\mu$. The heuristic interpretation will be that a constant mutation rate $\mu$ erodes all families at the same rate. The mutated lineages are unique and each one only takes up an infinitesimal fraction of the whole population until they possibly increase their size to a positive fraction of the population by a jump. They could also experience yet another mutation but that does not matter in this setting since we are not interested in the actual genotypes; all that matters is that they differ. In the case with finite intensity of births of new litters, the jump mechanism will be the same as in \eqref{rho_jump}, but for $t$ between two consecutive jump times, say $\sigma$ and $\tau$, we will have
\begin{equation}\label{rho_erodes}
\rho_t=e^{-\mu(t-\sigma)}\rho_{\sigma}+(1-e^{-\mu(t-\sigma)})\lambda,
\end{equation}
where $\lambda$ is the Lebesgue measure on $[0,1]$.

To make this rigorous, we will proceed in several steps. We will first study the litters without any genealogical relationships. Since a family consists of litters, claiming that it erodes at a constant rate, implies that its litters also must do so at the same rate. We will describe the composition of a population consisting of eroding litters and ``mutants'', or singletons, with a probability measure on $[0,\infty)$. The process describing the evolution of the population, still disregarding possible family ties between litters, will then be shown to converge to a stationary distribution. After that, we will impose a genealogy on the litters, meaning a partial order describing who is a descendant of whom. This will enable us to define $\rho$ as we want. Finally, Theorem \ref{main} validates our construction by stating that a sample from this population would have the same sampling distribution as from a $\Lambda$-coalescent with mutations.

We will for the rest of this section assume that $\Lambda(0)=\Lambda(1)=0$ and $\int_{(0,1)}x^{-1}\Lambda(dx)<\infty$. We let $\nu(dx):=\Lambda(dx)/x^2$, as in Section \ref{section:paintboxes_and_bridges}, and thus $\int_{(0,1)}x\nu(dx)<\infty$. Let $\{(\tau_i,X_i,U_i)\}_{i\in\N}$ be a Poisson process on $\R\times(0,1)\times(0,1)$ with intensity measure $dt\otimes\nu(dx)\otimes du$, whose points denote the times of birth, and the sizes of the litters born in the population, and auxiliary random variables for each litter to be used later. The litters are indexed in decreasing order of size, and in case of ties in decreasing order of the auxiliary random variables. Thus $i<j$ need not imply $\tau_i<\tau_j$. The dynamics in \eqref{rho_jump}, when there are no mutations, imply that the fraction of the population at time $t$ that belongs to a litter $i$, born at time $\tau_i\leq t$, with original size $X_i$, is given by
$$
X_i\prod_{j:\tau_i<\tau_j\leq t}(1-X_j),
$$
since the size of the litter must be scaled down at the birth of each subsequent litter. If the size of each litter, and thus also the size of each family, is furthermore eroded with a constant rate $\mu$, the size at time $t$ of litter $i$ becomes
$$
X_i(t):= X_ie^{-\mu(t-\tau_i)}\prod_{j:\tau_i<\tau_j\leq t}(1-X_j)\bbone(\tau_i \leq t),
$$
where $\bbone(\cdot)$ is the indicator function. We can describe the sizes and the ages of the litters at time $t$ with the random probability measure $\tilde{\rho}_t$ on $\R_+$, defined by its distribution function
\begin{equation*}
F(s:{\tilde{\rho}_t}):=1-e^{-\mu s}\prod_{i:t-s\leq \tau_i\leq t}(1-X_i),
\end{equation*}
for $s\geq 0$ and $F(s:{\tilde{\rho}_t}):=0$ for $s<0$. The atoms of $\tilde{\rho}_t$ now have sizes $X_i(t)$ and positions $t-\tau_i$, provided $\tau_i\leq t$, corresponding to the current sizes and ages at time $t$ of the litters. By the homogeneity of the Poisson process, $\tilde{\rho}=(\tilde{\rho}_t)_{t\in\R}$ is a stationary process. Since the process depends on \emph{all} litters born before time $t$, it  does not describe the composition of the population into litters if we want the process to start at time 0 with \emph{no} litters. In that case we must use a cut-off, so that there are no litters older than $t$ at time $t$. This is described by the process $\tilde{\rho}'=(\tilde{\rho}'_t)_{t\geq 0}$ of random probability measures on $\R_+$, with distribution functions
$$
F(s:\tilde{\rho}'_t):=1-e^{-\mu s}\prod_{i:\max(t-s,0)\leq \tau_i\leq t}(1-X_i),
$$
for $s\geq 0$ and $F(s:{\tilde{\rho}'_t}):=0$ for $s< 0$, so that $\tilde{\rho}'_0(ds)=\mu e^{-\mu s}ds$. This process is not stationary, but it converges in distribution to $\tilde{\rho}_0$.
\begin{lemma}\label{tilde_rho_conv}
$\tilde{\rho}'_t\overset{d}{\to}\tilde{\rho}_0$, as $t\to\infty$.
\end{lemma}
{\bf Proof} Define $\tilde{\rho}''_t,t\geq 0$, by its distribution function
\begin{equation*}
F(s:{\tilde{\rho}''_t}):=1-e^{-\mu s}\prod_{i:\max(-s,-t)\leq \tau_i\leq 0}(1-X_i),
\end{equation*}
for $s\geq 0$ and $F(s:{\tilde{\rho}''_t}):=0$ for $s<0$. By the homogeneity of the Poisson process, $F(s:{\tilde{\rho}'_t})\overset{d}{=}F(s:{\tilde{\rho}''_t})$, or, equivalently, $\tilde{\rho}'_t\overset{d}{=}\tilde{\rho}''_t$, for all fixed $t\geq 0$. Note that $F(s:{\tilde{\rho}_0})=F(s:{\tilde{\rho}''_t})$ for $s\leq t$. For $s>t$, we have
$$
0\leq F(s:{\tilde{\rho}_0})-F(s:{\tilde{\rho}''_t})\leq e^{-\mu s} < e^{-\mu t}\to 0,
$$
as $t\to\infty$. Thus, the distance between $\tilde{\rho}''_t$ and $\tilde{\rho}_0$ in the total variation metric, $\sup_{A\in\mathcal{B}(\R)}|\tilde{\rho}''_t(A)-\tilde{\rho}_0(A)|<e^{-\mu t}\to 0$, as $t\to\infty$. (Here $\mathcal{B}(\R)$ are the Borel sets of $\R$.) This implies $\tilde{\rho}'_t\overset{d}{\to}\tilde{\rho}_0$.\hfill$\Box$

Thus $\tilde{\rho}$ and $\tilde{\rho}'$ have the same limiting distribution, and we choose to work with the former process, since it is stationary.
\begin{theorem}\label{litter_comp}
The composition of a sample from $\tilde{\rho}_0$, according to litters of increasing age, is regenerative with characteristics $(\mu,\nu)$.
\end{theorem}
{\bf Proof} 
By construction, $F(s:\tilde{\rho}_0)$ is  a multiplicative subordinator with characteristics $(\mu,\nu)$, and the order of the parts of the regenerative composition obtained by sampling from the closure of its range corresponds to increasing age of the litters.\hfill$\Box$

In the light of Proposition \ref{mohle_neg}, the theorem might be a bit surprising. What we really want is not the composition into litters, but the partition into families, so we must somehow collect different litters into families. This will destroy the regenerative property of the composition into litters.

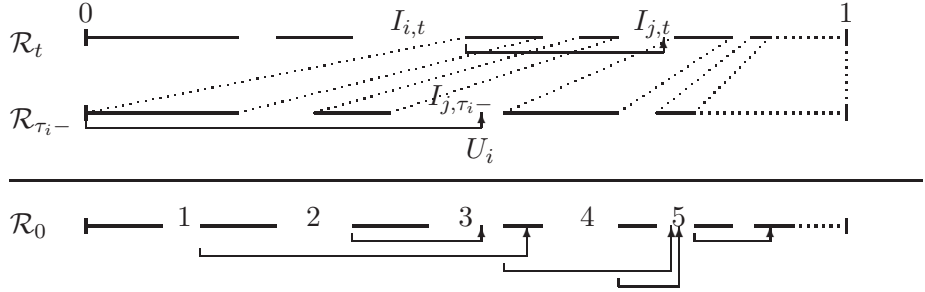
\begin{figure}
\centering
\setlength{\unitlength}{1cm}
\begin{picture}(12,4)(-1,-1.3)
\thicklines
\drawline(0,1.9)(0,2.1)
\drawline(10,1.9)(10,2.1)
\put(-1,1.8){$\mathcal{R}_t$}
\put(-0.1,2.2){0}
\put(9.9,2.2){1}
\drawline(0,2)(2,2)
\drawline(2.5,2)(3.5,2)
\drawline(5,2)(6,2)
\drawline(6.5,2)(7,2)
\drawline(7.75,2)(8.5,2)
\drawline(8.75,2)(9,2)
\dottedline{0.1}(9,2)(10,2)
\thinlines

\put(4,2.1){$I_{i,t}$}
\put(7.2,2.1){$I_{j,t}$}
\drawline(5,1.8)(5,1.9)
\drawline(5,1.8)(7.6,1.8)
\put(7.6,1.8){\vector(0,1){0.2}}
\put(4.5,1.1){$I_{j,\tau_i-}$}
\drawline(0,0.8)(5.2,0.8)
\drawline(0,0.8)(0,0.9)
\put(5.2,0.8){\vector(0,1){0.2}}
\put(5,0.4){$U_i$}

\dottedline{0.1}(5,2)(0,1)
\dottedline{0.1}(6,2)(2,1)
\dottedline{0.1}(6.5,2)(3,1)
\dottedline{0.1}(7,2)(4,1)
\dottedline{0.1}(7.75,2)(5.5,1)
\dottedline{0.1}(8.5,2)(7,1)
\dottedline{0.1}(8.75,2)(7.5,1)
\dottedline{0.1}(9,2)(8,1)
\dottedline{0.1}(10,1.9)(10,1.1)

\thicklines
\drawline(0,0.9)(0,1.1)
\drawline(10,0.9)(10,1.1)
\put(-1,0.8){$\mathcal{R}_{\tau_i-}$}
\drawline(0,1)(2,1)
\drawline(3,1)(4,1)
\drawline(5.5,1)(7,1)
\drawline(7.5,1)(8,1)
\dottedline{0.1}(8,1)(10,1)
\thinlines
\drawline(-1,0.1)(11,0.1)
\thicklines
\drawline(0,-0.6)(0,-0.4)
\drawline(10,-0.6)(10,-0.4)
\put(-1,-0.6){$\mathcal{R}_0$}
\drawline(0,-0.5)(1,-0.5)
\drawline(1.5,-0.5)(2.5,-0.5)
\drawline(3.5,-0.5)(4.5,-0.5)
\drawline(5.5,-0.5)(6,-0.5)
\drawline(7,-0.5)(7.5,-0.5)
\drawline(8,-0.5)(8.5,-0.5)
\drawline(8.8,-0.5)(9.3,-0.5)
\dottedline{0.1}(9.3,-0.5)(9.9,-0.5)
\thinlines
\put(1.2,-0.5){1}
\put(2.9,-0.5){2}
\put(4.9,-0.5){3}
\put(6.5,-0.5){4}
\put(7.7,-0.5){5}
\drawline(1.5,-0.9)(5.8,-0.9)
\drawline(1.5,-0.9)(1.5,-0.8)
\put(5.8,-0.9){\vector(0,1){0.4}}
\drawline(3.5,-0.7)(5.2,-0.7)
\drawline(3.5,-0.7)(3.5,-0.6)
\put(5.2,-0.7){\vector(0,1){0.2}}
\drawline(5.5,-1.1)(7.7,-1.1)
\drawline(5.5,-1.1)(5.5,-1)
\put(7.7,-1.1){\vector(0,1){0.6}}
\drawline(7,-1.3)(7.8,-1.3)
\drawline(7,-1.3)(7,-1.2)
\put(7.8,-1.3){\vector(0,1){0.8}}
\drawline(8,-0.7)(9,-0.7)
\drawline(8,-0.7)(8,-0.6)
\put(9,-0.7){\vector(0,1){0.2}}
\end{picture}
\caption{\textsc{Top} An illustration of $j\prec'i$. The arrow indicates how $U_i\in I_{j,\tau_i-}$.\newline
\textsc{Bottom} An example of how the litters $1,\dots,5$ (an arbitrary enumeration) can be related to each other. Here $5\prec'4$, $5\prec'3\prec'2$ and both $1$ and $5$ are roots.}\label{Rt}
\end{figure}

We will now define how the litters are related to each other. We do this by sampling from $\tilde{\rho}$. Let $\mathcal{R}_t$ be the closed range of $F(s:\tilde{\rho}_t)$. The complement of $\mathcal{R}_t$ in $[0,1]$ is a union of disjoint open intervals, $\cup_i I_{i,t}$, with
$$
I_{i,t}:=(F((t-\tau_i)-:\tilde{\rho}_t),F(t-\tau_i:\tilde{\rho}_t)),
$$
so that interval $I_{i,t}$ corresponds to litter $i$. Note that litters $i$ with $\tau_i> t$, i.e.\ litters not yet born at time $t$, have $I_{i,t}=\emptyset$. We also have $\mathcal{R}_t^{(u)}=\mathcal{R}_{\tau_i-}$ if $u\in I_{i,t}$, see Figure \ref{Rt}.

We say that litter $i$ originates from litter $j$ if $U_i\in I_{j,\tau_i-}$, and in that case we write $j\prec' i$, see Figure \ref{Rt}. There is for each $i$ at most one $j$ such that $j\prec'i$. There is no such $j$ if $U_i\in\mathcal{R}_{\tau_i-}$. Let $\mathcal{I}_0:=\{i:\nexists j,j\prec' i\}$. This is the set of litters which are not descendants from any other litter, but descendants from singletons, thus their genotypes are unique at their times of births. We call these litters \emph{roots}. We define $\prec$ by $j\prec i$ if  there exist $k_1,\dots,k_n$ such that $j\prec'k_1\prec'\cdots\prec'k_n= i$. The sequence $k_1,\dots,k_n$ is then unique. Furthermore, we set $j\preceq i$ if $j\prec i$ or $i=j$. There can be at most one root $j$ for each $i$ such that $j\preceq i$, and in that case we write $\alpha(i)=j$, and say that $j$ is the root of $i$. What is not immediately obvious is that each $i\in\N$ has a root (almost surely).
\begin{lemma}\label{grounded}
Each $i\in\N$ has a root almost surely.
\end{lemma}
\noindent{\bf Proof} Define recursively $\mathcal{I}_n:=\{i:\exists j\in\mathcal{I}_{n-1},j\prec'i\}$, for $n\geq 1$, and let the \emph{height} of a fixed litter $i$ be defined by $H_i:=n$ if $i\in\mathcal{I}_n$ and $H_i=\infty$ if $\nexists n\in\N:i\in\mathcal{I}_n$. We need to show that $H_i$ is finite almost surely. 

The height $H_i$ is a function of $\{(\tau_k,X_k,U_k)\}_{k:\tau_k<\tau_i}$. The event $\{i\prec' l\}$ is likewise measurable with respect to $\{(\tau_k,X_k,U_k)\}_{k:\tau_i\leq\tau_k<\tau_l}$. Thus the events $\{H_i=h\}$ and $\{i\prec'l\}$ are independent for all $i,l:\tau_i<\tau_l$, and $h\in\N\cup\{\infty\}$.

Let $\bar{p}_h:=P(H_i\geq h)$ for $h\geq 0$. Obviously $\bar{p}_0=1$. Assume $h\geq 1$. Note that $\{H_i\geq 1\}=\{U_i\notin\mathcal{R}_{\tau_i-}\}=\cup_{j:\tau_j<\tau_i}\{j\prec'i\}$, so
\begin{align*}
\bar{p}_h &= \sum_{j:\tau_j<\tau_i} P(H_i\geq h,j\prec'i)= \sum_{j:\tau_j<\tau_i} P(H_j\geq h-1,j\prec'i) \\
&= \sum_{j:\tau_j<\tau_i} P(H_j\geq h-1)P(j\prec'i) = \sum_{j:\tau_j<\tau_i} \bar{p}_{h-1}P(j\prec'i) \\
&=\bar{p}_{h-1}P(U_i\notin\mathcal{R}_{\tau_i-})=\bar{p}_{h-1}(1-q(1:1)').
\end{align*}
In the last equality we used \eqref{qn1prim} with $n=1$. Thus $\bar{p}_h=(1-q(1:1)')^h$ and therefore $H_i$ has a geometric distribution with parameter $q(1:1)'$ and is finite almost surely.\hfill$\Box$

By our interpretation of the relation $\preceq$ as a genealogical relation, we should let all litters with the same root have the same genotype. We define $R_i$, the genotype of litter $i$, by $R_i:=U_{\alpha(i)}$. Now we can finally define $\rho=(\rho_t)_{t\in\R}$.
$$
\rho_t:=\sum_iX_i(t)\delta_{R_i}+\Big(1-\sum_iX_i(t)\Big)\lambda.
$$
Note that this is a stationary version, and $\rho_0\not\equiv\lambda$. In the finite intensity case, $\rho$ behaves as in $\eqref{rho_erodes}$ between jumps, just as we wanted, and at the time of a jump, the new litter chooses its genotype from the population at the moment before the jump, just as in \eqref{rho_jump}.

At a fixed time $t$, $\rho_t$ represents the population in the sense that a sample from the population will have a partition with distribution as given by $\sim_{\rho_t}$. An i.i.d.\ sample $(r_i)_{i\in[n]}$ from a realization of $\rho_t$ can be interpreted as the genotypes of individuals $i=1,\dots,n$ in a sample from the population a time $t$. The value of an $r$ with distribution $\rho_t$ can either be one of the $R_1,\dots$, or, with probability $1-\sum_i X_i(t)$, it is uniformly drawn from $[0,1]$.

The justification for the construction is given by Theorem \ref{main}.
\begin{theorem}\label{main}
The partition of a sample from $\rho_0$ according to families has the same distribution as the partition according to genotypes of a sample from a $\Lambda$-coalescent with mutations, i.e.\ its distribution is given by the recursion \eqref{Mohles_recursion}.
\end{theorem}
\noindent{\bf Proof} We assume the sample size $n\geq 2$ and that the sample is created by first sampling from $\mathcal{R}_0$ with the i.i.d. uniform random variables $(V_i)_{i\in[n]}$, and then collecting the litters into families. Note that $(V_{in})_{i=j\dots n}$, when disregarding their order, are i.i.d.\ $U(v,1)$, given $V_{jn}>v$. We will use the notation from Section \ref{section:regenerative_composition_structures}. Consider the realization $V_{1n}=v$. Three possibilities exist. 
\begin{enumerate}
\item\label{ettan} $v\in\mathcal{R}_0$. This happens with probability $q(n:1)'$. Then 1 is a singleton and is thus in a family of its own.
\item\label{tvaan} $V_{1n}\notin\mathcal{R}_0$ and $[V_{1n},V_{2n}]\cap\mathcal{R}_0\neq\emptyset$. This happens with probability $q(n:1)''$. Then 1 is in a litter of its own, say litter $k$, and what family litter $k$ belongs to is determined by a uniform random variable $U_k$ on $\mathcal{R}^{(v)}_0=\mathcal{R}_{\tau_k-}$.
\item\label{trean} $[V_{1n},V_{mn}]\cap\mathcal{R}_0=\emptyset$ and either $m=n$, or $2\leq m< n$ and $[V_{mn},V_{m+1,n}]\cap\mathcal{R}_0\neq\emptyset$. Then $m$ individuals belong to the same litter, say litter $k$. This happens with probability $q(n:m)$. What family this litter belongs to is determined by a uniform random variable on $\mathcal{R}^{(v)}_0=\mathcal{R}_{\tau_k-}$.
\end{enumerate}
In case \ref{ettan}., we immediately find that the first part of our composition has size 1. The distribution of the rest of the sample is determined by $(V_{in})_{i=2\dots n}$ and $\mathcal{R}^{(v)}_0$. By the regenerative property, the distribution of the rest of the sample will be the same as sampled with $(V_i)_{i\in[n-1]}$ from $\mathcal{R}_0$.

In case \ref{tvaan}., the lineages of the sample can be represented by $U_k$ and $(V_{in})_{i=2\dots n}$ and their partition is obtained by sampling from $\mathcal{R}^{(v)}_0=\mathcal{R}_{\tau_k-}$, which by the regenerative property yields the same result in distribution as sampling with $(V_i)_{i\in[n]}$ from $\mathcal{R}_0$.

In the third case, we know that the lineages represented by $(V_{in})_{i\in[m]}$ have coalesced since they originate from a common litter, say litter $k$, but we do not know to which family they belong. This is determined by the realization of $U_k$ relative to $\mathcal{R}^{(v)}_0$, which, if $2\leq m< n$, together with the realization of $(V_{in})_{i=m+1,\dots,n}$ determines the further coalescing of lineages. As in case \ref{tvaan}., the distribution will be the same as if we sample with $(V_i)_{i\in[n-m+1]}$ from $\mathcal{R}_0$.

The argument is now similar to the one in Möhle \cite{Mohle:2005a} and Dong et al.\ \cite{DongEtAl:2006}, with the main difference that our case \ref{tvaan}.\ above does not add any information about the final partition, whereas they only have either mutations/freezing (our case \ref{ettan}.) or collisions (our case \ref{trean}.) happening at each stage. We thus get the recursion
\begin{align*}
q(\mathbf{a}) &= q(n:1)' q(\mathbf{a}-\mathbf{e}_1) + q(n:1)''q(\mathbf{a})\\
&\qquad+\sum_{m=2}^n q(n:m)\sum_{j=1}^{n-m+1}\frac{j(a_j+1)}{n-m+1}q(\mathbf{a}+\mathbf{e}_j-\mathbf{e}_{j+m-1}), \\
\intertext{or equivalently, by \eqref{qnm}, \eqref{qn1prim}, \eqref{qn1bis} and $\Phi(n:m)=\binom{n}{m}\lambda_{n,m}$ for $2\leq  m\leq n$,}
\Phi(n)q(\mathbf{a}) &= \mu n q(\mathbf{a}-\mathbf{e}_1) + n\int_0^1x(1-x)^{n-1}\nu(dx)q(\mathbf{a})\\
&\qquad + \sum_{m=2}^n \binom{n}{m}\lambda_{n,m}\sum_{j=1}^{n-m+1}\frac{j(a_j+1)}{n-m+1}q(\mathbf{a}+\mathbf{e}_j-\mathbf{e}_{j+m-1}),
\end{align*}
and since $\Phi(n)-n\int_0^1x(1-x)^{n-1}\nu(dx)=\mu n + \sum_{m=2}^n\binom{n}{m}\lambda_{n,m}=\mu n+\lambda_n$, we arrive at \eqref{Mohles_recursion}, and the proof is complete.\hfill$\Box$

\begin{figure}
\centering
\setlength{\unitlength}{1cm}
\begin{picture}(11,1.8)(-1,-0.8)
\thicklines
\drawline(0,-0.1)(0,0.1)
\drawline(10,-0.1)(10,0.1)
\put(-1,-0.2){$\mathcal{R}_0$}
\put(-0.1,-0.5){0}
\put(9.9,-0.5){1}
\drawline(0,0)(1,0)
\drawline(1.5,0)(2.5,0)
\drawline(3.5,0)(4.5,0)
\drawline(5.5,0)(6,0)
\drawline(7,0)(7.5,0)
\drawline(8,0)(8.5,0)
\drawline(8.8,0)(9.3,0)
\dottedline{0.1}(9.3,0)(9.9,0)
\thinlines
\drawline(1.5,-0.4)(5.8,-0.4)
\drawline(1.5,-0.4)(1.5,-0.3)
\put(5.8,-0.4){\vector(0,1){0.4}}
\drawline(3.5,-0.2)(5.2,-0.2)
\drawline(3.5,-0.2)(3.5,-0.1)
\put(5.2,-0.2){\vector(0,1){0.2}}
\drawline(5.5,-0.6)(7.7,-0.6)
\drawline(5.5,-0.6)(5.5,-0.5)
\put(7.7,-0.6){\vector(0,1){0.6}}
\drawline(7,-0.8)(7.8,-0.8)
\drawline(7,-0.8)(7,-0.7)
\put(7.8,-0.8){\vector(0,1){0.8}}
\drawline(8,-0.2)(9,-0.2)
\drawline(8,-0.2)(8,-0.1)
\put(9,-0.2){\vector(0,1){0.2}}
\put(1.2,0.2){\vector(0,-1){0.2}}
\put(1.1,0.3){$V_4$}
\put(2,0.2){\vector(0,-1){0.2}}
\put(1.9,0.3){$V_7$}
\put(2.7,0.2){\vector(0,-1){0.2}}
\put(2.6,0.3){$V_2$}
\put(3.3,0.2){\vector(0,-1){0.2}}
\put(3.2,0.3){$V_5$}
\put(6.15,0.2){\vector(0,-1){0.2}}
\put(6.00,0.3){$V_1$}
\put(6.55,0.2){\vector(0,-1){0.2}}
\put(6.40,0.3){$V_3$}
\put(6.9,0.2){\vector(0,-1){0.2}}
\put(6.85,0.3){$V_6$}
\end{picture}
\caption{Illustration of sampling with $(V_i)_{i\in[7]}$ from the regenerative set that corresponds to $\tilde{\rho}_0$. The arrows indicate how the litters are related to each other. Compare with Figures \ref{Rz} and \ref{Rt}.}\label{Example}
\end{figure}
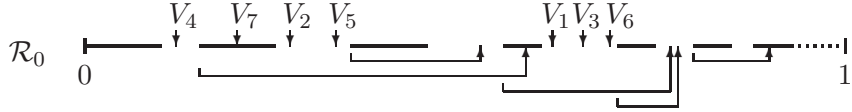
We illustrate the procedure of the proof with Figure \ref{Example}. The procedure amounts to moving from left to right and note how the arrows hit $\mathcal{R}_0$ or its complement. First, $V_4$ is alone in its interval, corresponding to case \ref{tvaan}.\ in the proof. At this point we cannot say anything about the final partition since that litter may be related to the other individuals in our sample. Next, $V_7$ hits $\mathcal{R}_0$ so that it is in a family of its own (case \ref{ettan}.). The next event is that both $V_2$ and $V_5$ fall in the same interval, corresponding to a merger of their lineages and case \ref{trean}.\ in the proof. After that, we have a case \ref{tvaan}.\ for the lineage of 2 and 5. Next, we find that the litter of individual 4 is a root. Then lineages 1, 3 and 6 coalesce. The penultimate event is that the litters of lineages 2 and 5, and 1, 3 and 6, are related to the same litter, and thus these lineages coalesce. The final event is finding that this litter also is a root. Thus the partition is $\big\{\{1,2,3,5,6\},\{4\},\{7\}\big\}$, just as the example of Figure \ref{FamilyTree}. The order of the collisions and mutations is also the same as in that example.

\begin{remark}\label{concl} Our construction of $\rho$ requires $\nu$ to be a measure on $(0,1)$ with $\int_{(0,1)}x\nu(dx)<\infty$. This excludes a large class of $\Lambda$-coalescents. The moment condition is necessary when we want to construct the multiplicative subordinator $F(s:\tilde{\rho}_0)$ (whose properties we use repeatedly) from the point process $\{(\tau_i,X_i)\}_{i\in\N}$. Nevertheless, it might be possible to obtain a convergence result analogous to the one of Proposition \ref{BertoinLeGall}, but we have not been able to do so.
\end{remark}

\noindent\textbf{Acknowledgment} I thank the two referees for their thorough reading and detailed comments that certainly helped to improve the presentation and rigor of this paper.

\end{document}